\documentclass[a4paper,11pt]{amsart}
\usepackage{amssymb}

\theoremstyle{plain}
\newtheorem{theorem}{Theorem}[section]

\newtheorem{proposition}[theorem]{Proposition}

\theoremstyle{definition}

\theoremstyle{remark}
\newtheorem{remark}[theorem]{Remark}
\newcommand{\kc}{{\mathcal C}}
\newcommand{\kd}{{\mathcal D}}

\newcommand{\ko}{{\mathcal O}}
\newcommand{\kp}{{\mathcal P}}

\newcommand{\kr}{{\mathcal R}}
\newcommand{\IA}{{\mathbb A}}

\newcommand{\IC}{{\mathbb C}}
\newcommand{\IH}{{\mathbb H}}

\newcommand{\IN}{{\mathbb N}}
\newcommand{\IP}{{\mathbb P}}
\newcommand{\IQ}{{\mathbb Q}}

\newcommand{\IZ}{{\mathbb Z}}

\newcommand{\gothq}{{\mathfrak q}}
\newcommand{\id}{{\rm id}}
\newcommand{\tensor}{\otimes}

\newcommand{\Hilb}{{\rm Hilb}}
\newcommand{\lra}{\longrightarrow}
\newcommand{\xra}{\xrightarrow}
\newcommand{\Hilbn}{\Hilb^n(\IA^2_\IC)}

\newcommand{\ie}{i.e.\ }
\DeclareMathOperator{\Ind}{Ind}
\newcommand{\eps}{\varepsilon}
\DeclareMathOperator{\sgn}{sgn}

\begin{document}
\title[Cup Product]{Symmetric Groups and the Cup Product on the
Cohomology of Hilbert Schemes}
\author{Manfred Lehn \and Christoph Sorger}

\begin{abstract}
Let $\kc(S_n)$ be the $\IZ$-module of integer valued class functions
on the symmetric group $S_n$. We introduce a graded version of the
convolution product on $\kc(S_n)$ and show that there is a degree
preserving ring isomorphism $\kc(S_n)\lra H^*(\Hilbn;\IZ)$ to the
cohomology of the Hilbert scheme of points in the complex affine
plane.
\end{abstract}
\address{
Manfred Lehn\\
Mathematisches Institut
der Universit\"at zu K\"oln\\
Weyertal 86-90\\
D-50931 K\"oln, Germany}
\email{manfred.lehn@math.uni-koeln.de}
\address{
Christoph Sorger\\
Math\'{e}matiques (UMR 6629 du CNRS)\\
Universit\'{e} de Nantes\\
2, Rue de la Houssini\`{e}re\\
BP 92208\\
F-44322 Nantes Cedex 03, France}
\email{christoph.sorger@math.univ-nantes.fr}


\maketitle

\section{Introduction}
In this paper we relate a geometric and a group theoretic incarnation
of the bosonic Fock space $\kp=\IQ[p_1,p_2,p_3,\ldots]$. On the
geometric side this is the direct sum
$$\IH=\bigoplus_{n\geq 0}H^*(\Hilbn;\IQ)$$
of the rational cohomology of the Hilbert schemes of generalized
$n$-tuples in the complex affine plane, and on the group theoretic
side the direct sum
$$\kc=\bigoplus_{n\geq 0}\kc(S_n)\tensor_\IZ\IQ$$
of the spaces of class functions on the symmetric groups $S_n$. In
addition to their vertex algebra structures, both $\IH$ and $\kc$
carry natural ring structures on each component of fixed conformal
weight: for the cohomology of the Hilbert schemes this is the
ordinary topological cup product; for $\kc$ it is a certain
combinatorial cup product to be defined below (equation
(\ref{eq:DefOfCup})) and not to be confused with the usual product
on the representation ring arising from tensor product of
representations. With respect to these  products we can state our
main theorem:

\begin{theorem}\label{th:Main}---
The composite isomorphism of vertex algebras
$$\kc\stackrel{\Phi}\lra \kp\stackrel{\Psi}\lra \IH$$
induces for each conformal weight $n\in \IN_0$ an isomorphism of
graded rings
$$\kc(S_n)\lra H^*(\Hilbn;\IZ).$$

\end{theorem}

That the composition $\Psi\Phi$ is an isomorphism of vertex algebras
is by now well known (see \cite{Nakajima} for $\Psi$ and \cite{FJW}
for $\Phi$). The emphasis of the theorem is on the multiplicativity
of this map. Based on the geometric analysis carried out in
\cite{Lehn}, the proof is purely algebraic. Our starting point was
the observation of Frenkel and Wang \cite{Frenkel-Wang} that
Goulden's differential operator $\Delta$ \cite{Goulden} is closely
related to the operator $\partial$ of \cite{Lehn}.

The first named author gratefully acknowledges the support and the
hospitality of the University of Nantes, where this paper was
written.

{\small Note added in proof: E. Vasserot [arXiv:math.AG/0009127] has
independently obtained a similar result using other methods.}

\section{The cup product on the group ring $\IZ[S_n]$}

Let $\kc(S_n)$ denote the set of integer valued class functions on
the symmetric group $S_n$, \ie the set of functions $S_n\to\IZ$ which
are constant on conjugacy classes. Identifying a function $f$ with
the linear combination $\sum_{\pi\in S_n}f(\pi)\pi$, we may think of
$\kc(S_n)$ as a $\IZ$-submodule of the group ring $\IZ[S_n]$. As
such it inherits a product
$$(f*g)(\pi)= \sum_{\sigma\in S_n}f(\pi\sigma^{-1})g(\sigma),$$
called the convolution product.

\begin{remark}--- The character map
$\chi:\kr(S_n)\tensor_{\IZ}\IQ\lra\kc(S_n)\tensor_{\IZ}\IQ$ is a
$\IQ$-linear isomorphism from the rational representation ring to the
ring of rational class functions. The tensor product ring structure
on $\kr(S_n)$ is quite different from the convolution product
structure on $\kc(S_n)$, so that $\chi$ is {\em not} a ring
homomorphism. Even though we will use the identification of
$\kr(S_n)\tensor_{\IZ}\IQ$ and $\kc(S_n)\tensor_{\IZ}\IQ$ in the
definition of the vertex algebra structure on $\bigoplus_{n\geq0}
\kc(S_n)\tensor_{\IZ}\IQ$, we will never use the tensor product ring
structure in this paper.
\end{remark}

An integral basis of $\kc(S_n)$ is given by the characteristic
functions
$$\chi_\lambda=\sum_{\text{$\pi$ of type $\lambda$}}\pi,$$
where $\lambda$ is a partition of $n$ and $\pi$ runs through all
permutations with cycle type $\lambda$, \ie those having a disjoint
cycle decomposition with cycle lengths
$\lambda_1,\lambda_2,\ldots,\lambda_s$. For instance, the unit
element of the group ring $\IZ[S_n]$ -- and of $\kc(S_n)$ -- is
$\chi_{[1,1,\ldots,1]}$.

For any partition $\lambda$, let $\ell(\lambda)$ denote the length of
$\lambda$.
We introduce a gradation
$$\IZ[S_n]=\bigoplus_{d=0}^{n-1}\IZ[S_n](d)$$
as follows: a permutation $\pi$ has degree $\deg(\pi)=d$ if it can be
written as a product of $d$ transpositions but not less.
Equivalently, if $\pi$ is of cycle type $\lambda$, then $\deg(\pi)=n-\ell
(\lambda)$.

In particular, the maximal possible degree is indeed $n-1$. The
product in $\IZ[S_n]$ does not preserve this gradation, but it is
clearly compatible with the associated filtration
$$F^{d}\IZ[S_n]:=\bigoplus_{d'\leq d}\IZ[S_n](d'),$$
\ie it satisfies
$$F^{i}\IZ[S_n]*F^{j}\IZ[S_n]\subset F^{i+j}\IZ[S_n].$$
The induced product on
$\IZ[S_n]=gr^F\IZ[S_n]=\bigoplus_{d=0}^{n-1}F^d\IZ[S_n]/F^{d-1}\IZ[S_n]$
will be called {\em cup product} and denoted by $\cup$. Explicitly,
\begin{equation}\label{eq:DefOfCup}
\sigma\cup\pi=
\begin{cases}
\sigma*\pi&\text{ if $\deg(\sigma)+\deg(\pi)=\deg(\sigma\pi)$,}\\
0 &\text{ else.}
\end{cases}
\end{equation}
for $\pi,\sigma\in S_n$.

Clearly, the subring of class functions $\kc(S_n)\subset\IZ[S_n]$ is
generated by homogeneous elements and inherits from $\IZ[S_n]$
gradation, filtration, and, most importantly, the cup product.

Let  $\kc:=\bigoplus_{n\geq 0}\kc(S_n)\tensor_{\IZ}\IQ$. It is
bigraded by conformal weight $n$ and degree.

\section{The ring of symmetric functions}

Let $\kp=\IQ[p_1,p_2,p_3,\ldots]$ denote the polynomial ring in
countably infinitely many variables. It is endowed with a bigrading
by letting $p_m$ have conformal weight $m$ and cohomological degree
$m-1$. Let $\kp_n$ denote the component of conformal weight $n$, \ie
the subspace spanned by all monomials
$p_1^{\alpha_1}\cdot\ldots\cdot p_s^{\alpha_s}$ with
$\sum_ii\alpha_i=n$.

Define  linear maps $\Phi_n:\IQ[S_n]\lra \kp_n$ by sending a
permutation $\pi$ of cycle type
$\lambda=(\lambda_1\geq\lambda_2\geq\ldots\geq\lambda_s)$ to the
monomial $\frac{1}{n!}p_{\lambda_1}\cdot\ldots\cdot p_{\lambda_s}$.
Thus for any partition $\lambda=(1^{\alpha_1}2^{\alpha_2}\cdots)$ we
have
$$\Phi_n(\chi_\lambda)=\prod_i\frac{1}{\alpha_i!}
\left(\frac{p_i}{i}\right)^{\alpha_i}.$$
In particular, there is an
isomorphism of bigraded vector spaces
$$
\Phi: \kc \lra \kp.
$$
Moreover, multiplication by $p_m$ in $\kp$ corresponds to linear
operators $r_m$ in $\kc$ which are given as follows (see \cite{FJW}
and the references therein): let $\Ind$ denote induction of class
functions. Then $r_m$ is the map
$$
r_m:\kc(S_n)\xra{\id\tensor m\chi_{(m)}}
\kc(S_n)\tensor\kc(S_m)=\kc(S_n\times S_m) \xra{\Ind}\kc(S_{n+m}).
$$
The map $r_1$ is in fact very easy to describe: let $\iota:\IQ[S_n]
\to \IQ[S_{n+1}]$ be induced from the standard inclusion $S_n\to S_{n+1}$.
Then $r_1$ extends to a map
$$\IQ[S_n]\to\IQ[S_{n+1}],\quad \pi\mapsto \frac{1}{n!}\sum_{t\in
S_{n+1}}t\iota(\pi)
t^{-1}.$$

In \cite{Goulden}, I.P.~Goulden introduces the following differential
operator on $\kp$:
$$
\Delta  :=\Delta^\prime+\Delta^{\prime\prime} := \frac{1}{2}\sum_{i,j}
ijp_{i+j}\frac{\partial}{\partial p_i} \frac{\partial}{\partial p_j}
+\frac{1}{2} \sum_{i,j} (i+j)p_{i}p_{j}\frac{\partial}{\partial
p_{i+j}}
$$
and proves

\begin{proposition}[Goulden]--- Let $\tau_n\in\kc(S_n)$ denote the sum of
all transpositions in $S_n$. Then
$$\Phi(\tau_n*y)=\Delta(\Phi(y))$$
for all $y\in\IQ[S_n]$.
\end{proposition}

\begin{proof}
(\cite{Goulden}, Proposition 3.1)
\end{proof}

Note that $\Delta^{\prime}$ is of bidegree $(0,1)$ and that
$\Delta^{\prime\prime}$ is of bidegree $(0,-1)$. Since $\tau_n$ is of
degree one, for the cup product introduced above Goulden's
proposition reads as follows:
\begin{equation}\label{eq:GouldenCup}
\Phi(\tau_n\cup y)=\Delta^\prime(\Phi(y))
\end{equation}

\section{The Hilbert scheme of points}

Consider the Hilbert scheme $\Hilbn$ of generalized $n$-tuples of
points on the affine plane. Let us recall some basic facts: $\Hilbn$ is a
quasi-pro\-jec\-tive manifold of dimension $2n$ (\cite{Fogarty}). The
closed subset
$$\Hilbn_O=\{\xi\in\Hilbn\,|\,{\mathrm{Supp}}(\xi)=O\in\IA^2_\IC\}$$
is a deformation retract of $\Hilbn$. This subvariety is
$(n-1)$-dimensional and irreducible (\cite{Briancon}) and has a cell
decomposition with $p(n,n-i)$ cells of dimension $i$, where $p(n,j)$
denotes the number of partitions of $n$ into $j$ parts
(\cite{E-S2}). In particular, the odd dimensional cohomology
vanishes, there is no torsion and
$H^{2i}(\Hilbn;\IZ)=\IZ^{p(n,n-i)}$. In order not to worry
constantly about a factor of 2, we agree to give
$H^{2i}(\Hilbn;\IZ)$  degree $i$.

Consider the bigraded vector space
$$\IH:=\bigoplus_{0\leq n}\bigoplus_{0\leq i<n} H^{2i}(\Hilbn;\IQ).$$
We refer to $i$ as the cohomological degree and to $n$ as the
conformal weight. Nakajima (\cite{Nakajima}) and Grojnowski
(\cite{Grojnowski}) used incidence varieties to construct linear
operators $\gothq_m:\IH\lra\IH$ and an isomorphism
$$\Psi:\kp\lra \IH$$
such that $\Psi(p_m\cdot y)=\gothq_m(\Psi(y))$.

Now let $\Xi_n\subset \Hilbn\times\IA^2_\IC$ denote the universal
subscheme para\-meter\-ized by $\Hilbn$ and consider the direct image
$pr_{1*}\ko_{\Xi_{n}}$ of its structure sheaf under the projection
to the first factor. It was shown by Ellingsrud and Str\o{}mme
\cite{E-S1} that the components of the total Chern class
$\gamma_n:=c(pr_{1*}\ko_{\Xi_{n}})$ generate $H^*(\Hilbn;\IZ)$ as a
ring. The relations between these generators are encoded in the
following identity:

\begin{theorem}\label{th:BigDiffOperator}---
Consider the following differential operator on $\kp$:
$$
\kd:={\rm Coeff}\left(t^0, \left(-\sum_{m>0}p_mt^m\right)
\exp\left(-\sum_{m>0}m\frac{\partial}{\partial p_m}t^{-m}\right)
\right).
$$
Then for any $y\in\kp$ one has
\begin{equation}
ch(pr_{1*}\ko_{\Xi_{n}})\cup \Psi(y)=\Psi(\kd(y)).
\end{equation}
\end{theorem}

\begin{proof}(\cite{Lehn}, Theorem 4.10).\end{proof}

The degree one part of the operator $\kd$ is equal to
$-\Delta^\prime$, hence
\begin{equation}\label{eq:PartialDelta}
c_1(pr_{1*}\ko_{\Xi_{n}})\cup\Psi(y)=-\Psi(\Delta^\prime(y)).
\end{equation}

In order to simplify notations we write $\partial(y):=
c_1(pr_{1*}\ko_{\Xi_{n}})\cup y$ for $y\in H^*(\Hilbn;\IQ)$. We can
combine the identities (\ref{eq:GouldenCup}) and
(\ref{eq:PartialDelta}) to obtain:
\begin{equation}\label{eq:TauIdentity}
-\partial\Psi\Phi(y))=\Psi(\Delta^\prime(\Phi(y))=\Psi\Phi(\tau_n\cup y)
\end{equation}

\begin{proposition}\label{pr:Induktion}--- For each $n$ we have
\begin{enumerate}
\item $\kc(S_n)\tensor_{\IZ}\IQ=\tau_n\cup (\kc(S_n)\tensor_{\IZ}\IQ)+
r_1(\kc(S_{n-1})\tensor_{\IZ}\IQ)$.
\item $\kp_n=\Delta^\prime(\kp_n)+ p_1\kp_{n-1}$.
\item $H^*(\Hilbn;\IQ)=\partial\,H^*(\Hilbn;\IQ)+\gothq_1\,H^*(\Hilb^{n-1}
(\IA^2_\IC);\IQ)$.
\end{enumerate}
\end{proposition}

\begin{proof} In view of equation (\ref{eq:TauIdentity}) and the
fact that $\Phi$ and $\Psi$ are isomorphisms the three assertions are
of course equivalent. Assertion (2) follows from the identities
$$[\Delta^\prime,p_1]=\sum_{j>0}jp_{j+1}\frac{\partial}{\partial p_j}$$
and
$$ad([\Delta^\prime,p_1])^{n-1}(p_1)=(n-1)!\,p_n$$
by an easy induction.
\end{proof}

\begin{theorem}\label{th:Lehn}--- For all $y\in H^*(\Hilbn;\IQ)$ one has
\begin{equation}\label{eq:qcommutator}
\gamma_{n+1}\cup \gothq_1(y)-\gothq_1(\gamma_n\cup y)=[\partial,\gothq_1]
(\gamma_n \cup y)
\end{equation}
\end{theorem}

\begin{proof} This is Theorem 4.2 of \cite{Lehn}.\end{proof}

\section{The alternating character}

Let $\eps_n\in\kc(S_n)$ denote the alternating character, \ie
$$\eps_n=\sum_{\pi\in S_n}\sgn(\pi)\pi.$$

\begin{proposition} The following identity holds for all
$y\in\kc(S_n)\tensor_\IZ\IQ$:
\begin{equation}\label{eq:epscommutator}
\eps_{n+1}\cup r_1(y)-r_1(\eps_n\cup y)=-
\tau_{n+1}\cup r_1(\eps_n\cup y)+r_1(\tau_n\cup\eps_n\cup y)
\end{equation}
\end{proposition}

\begin{proof} First note that we have the identities
$$\displaystyle\tau_{n+1}-\iota(\tau_n)=\sum_{i=1}^{n}(i\
n+1),\quad\text{and}\quad
\displaystyle\eps_{n+1}-\iota(\eps_n)=-\sum_{i=1}^{n}(i\
n+1)\cup\iota(\eps_n),$$
which together give
\begin{equation}\label{nochnegleichung}
\eps_{n+1}-\iota(\eps_n)=(-\tau_{n+1}+\iota(\tau_n))\cup\iota(\eps_n).
\end{equation}
Then
\begin{eqnarray*}
\lefteqn{n!\,\Big[\eps_{n+1}\cup r_1(y)-r_1(\eps_n\cup y)\Big]}\hspace{3em}\\
&=&\left[\eps_{n+1}\cup\sum_t t\iota(y)t^{-1}-\sum_tt\iota(\eps_n\cup y)t^{-1}\right]\\
&=&\sum_tt[\eps_{n+1}-\iota(\eps_n)]\iota(y)t^{-1}\quad
\text{, since $\eps_{n+1}$ is symmetric}\\
&=&\sum_tt[-\tau_{n+1}+\iota(\tau_n)]\cup\iota(\eps_n\cup y)t^{-1}
\quad\text{, by (\ref{nochnegleichung})}\\
&=&-\tau_{n+1}\sum_tt\iota(\eps_n\cup y)t^{-1}+
\sum_tt\iota(\tau_n\cup\eps_n\cup y)t^{-1}\\
&=&n!\,\Big[-\tau_{n+1}\cup r_1(\eps_n\cup y)+r_1(\tau_n\cup\eps_n\cup y)\Big]
\end{eqnarray*}
\end{proof}

\begin{proposition}--- The following identities hold:
\begin{equation}
\sum_{n\geq 0}\Phi(\eps_n)z^n=\exp\left(\sum_{m>0}(-1)^{m-1}
\frac{z^m}{m}p_m\right)=\sum_{n\geq 0}\Psi^{-1}(\gamma_n)z^n
\end{equation}
\end{proposition}

\begin{proof} The first equality can be found in \cite{Frenkel-Wang}.
The second equality is Theorem 4.6 in \cite{Lehn}.
\end{proof}

Thus under the isomorphism $\Psi\Phi:\kc\to \IH$ the alternating
character $\eps_n$ is mapped to the total Chern class $\gamma_n$ of
the tautological sheaf $pr_{1*}\ko_{\Xi_n}$!

\begin{proposition}\label{pr:Main}--- For all $y\in\kc(S_n)\tensor_{\IZ}\IQ$
the following identity holds:
\begin{equation}\label{eq:Main}
\Psi\Phi(\eps_n\cup y) =\gamma_n\cup\Psi\Phi(y)
\end{equation}
\end{proposition}

\begin{proof} Because of Proposition \ref{pr:Induktion} we may assume
that $y$ is of the form $\tau_n\cup x$ or $r_1(x)$. We will therefore
argue by induction on weight and degree and assume that the assertion
holds for all $x$ of either less degree or less weight than $y$. (The
assertion is certainly trivial for the vacuum, the -- up to a scalar
factor -- unique element of weight 0 and degree 0.)

In case $y=\tau_n\cup x$ it follows from equation
(\ref{eq:TauIdentity}) and induction that
\begin{eqnarray*}
\Psi\Phi(\eps_n\cup y)&=&\Psi\Phi(\eps_n\cup \tau_n\cup x)\\
&=&\partial(\Psi\Phi(\eps_n\cup x))\\
&=&\partial(\gamma_n\cup\Psi\Phi(x))\\
&=&\gamma_n\cup \partial(\Psi\Phi(x))\\
&=&\gamma_n\cup \Psi\Phi(\tau_n\cup x)\\
&=&\gamma_n\cup\Psi\Phi(y).
\end{eqnarray*}

In case $y=r_1(x)$ for some $x\in\kc(S_{n-1})$
we argue as follows:% equation (\ref{eq:epscommutator}):
\begin{eqnarray*}
\Psi\Phi(\eps_{n}\cup y)&=&\Psi\Phi(\eps_{n}\cup r_1(x))\\
&=&\Psi\Phi(r_1(\eps_{n-1}\cup x)-\tau_n\cup r_1(\eps_{n-1}\cup x)\\
&&\quad+r_1(\tau_{n-1}\cup\eps_{n-1}\cup x))
\quad\mbox{by equation (\ref{eq:epscommutator})}\\
&=&\gothq_1\Psi\Phi(\eps_{n-1}\cup x)
+\partial\gothq_1(\Psi\Phi(\eps_{n-1}\cup x))\\
&&\quad-\gothq_1\partial(\Psi\Phi(\eps_{n-1}\cup x))
\quad\mbox{by equation (\ref{eq:TauIdentity})}\\
&=&(\gothq_1+\partial\circ\gothq_1-\gothq_1\circ\partial)(\gamma_{n-1}\cup
\Psi\Phi(x))\\
&&\quad\mbox{by induction}\\
&=&\gamma_n\cup \gothq_1(\Psi\Phi(x))\quad\mbox{by equation
(\ref{eq:qcommutator})}\\
&=&\gamma_n\cup\Psi\Phi(r_1(x))\\
&=&\gamma_n\cup\Psi\Phi(y).
\end{eqnarray*}
\end{proof}

\section{Proof of Theorem \ref{th:Main}}

It follows from Proposition \ref{pr:Main} that the assertion of Theorem
\ref{th:Main} holds for rational coefficients: from \cite{E-S1} we know
that the Chern classes of $pr_{1*}\ko_{\Xi_n}$ generate the ring
$H^*(\Hilbn;\IZ)$. Hence Proposition \ref{pr:Main} implies that the
isomorphism
$$\Psi\Phi:\kc(S_n)\tensor_{\IZ}\IQ\lra H^*(\Hilbn;\IQ)$$
preserves the cup product and that the homogeneous components of the
alternating character $\eps_n\in\kc(S_n)$ generate
$\kc(S_n)\tensor_{\IZ}\IQ$. Thus in order to prove Theorem \ref{th:Main}
it suffices to see that this holds as well over the integers:

\begin{proposition}--- The homogeneous components of the alternating
character $\eps_n$ generate the ring $\kc(S_n)$ of integer valued
class functions with respect to the cup product.
\end{proposition}

\begin{remark}---
Of course, this implies the analogous statement for $\kc(S_n)$
equipped with the convolution product.
\end{remark}

\begin{proof} Let $\eps_n(i)$ denote the component of $\eps_n$ of
degree $i$. We must show that for each $d\geq0$ the elements
\begin{equation*}
  \eps_n^\lambda:=\prod_{i\geq1}\eps_n(i)^{\alpha_i},
\end{equation*}
where $\lambda=(1^{\alpha_1}2^{\alpha_2}\ldots)$ runs through all
partitions of $d$, form a set of generators of the $\IZ$-module
$\kc(S_n)(d)$.

\medbreak\noindent{\em Claim: The restriction map
$\kc(S_n)\xra{\rho}\kc(S_{n-1})$ is a surjective and degree
preserving ring homomorphism and maps $\eps_n$ to $\eps_{n-1}$. In
particular, we may assume that $n\geq2d$.}

Indeed, surjectivity and homogeneity are obvious, hence it is enough
to check that $\rho$ is multiplicative. We have
\begin{eqnarray*}
\rho(g\cup f)(\pi)
&=& \sum_{\sigma\in S_n}^{\sim}g(\pi\sigma^{-1})f(\sigma)\\
&=& \sum_{\sigma\in S_{n-1}}^{\sim}g(\pi\sigma^{-1})f(\sigma)+
\sum_{\sigma\in S_n\setminus S_{n-1}}^{\sim}g(\pi\sigma^{-1})f(\sigma),
\end{eqnarray*}
where $\displaystyle\sum^\sim$ means the sum over terms satisfying
the degree condition (\ref{eq:DefOfCup}).

The first term of the last line equals $(\rho(g)\cup\rho(f))(\pi)$,
and it suffices to show that no summand of the second term occurs.
Indeed, we may decompose any $\sigma\in S_n\setminus S_{n-1}$ as
$\sigma=(i\,n)\eta$ for some transposition $(i\,n)$ with
$i\in\{1,\ldots,n-1\}$ and $\eta\in S_{n-1}$. The degree matching
condition requires
\begin{equation}\label{eq:matchcondition}
\deg(\pi)=\deg(\pi\eta^{-1}(i\,n))+\deg((i\,n)\eta).
\end{equation}
But the right hand side equals
\begin{equation}
\deg(\pi\eta^{-1})+1+\deg(\eta)+1\geq\deg(\pi)+2,
\end{equation}
so that (\ref{eq:matchcondition}) is never fulfilled. Finally,
it is clear that $\rho(\eps_n)=\eps_{n-1}$, which proves
the claim. Note that $\rho$ is not multiplicative with respect to the
convolution product.

\medbreak Assume from now on that $n\geq2d$ and consider the
$\IZ$-module $M=\Phi(\kc(S_n)(d))\subset\kp$. It has a $\IZ$-basis
consisting of monomials
\begin{equation}\label{eq:monomials}
  p_n^{\lambda}:=\prod_{i\geq1}\frac{1}{\alpha^{\prime}_{i}!}
\left(\frac{p_i}{i}\right)^{\alpha^{\prime}_{i}}=\Phi(\chi_{\lambda^\prime})
\end{equation}
where as before $\lambda=(1^{\alpha_1}2^{\alpha_2}\ldots)$ is a
partition of $d$ and
$\lambda^{\prime}=(1^{\alpha^{\prime}_1}2^{\alpha^{\prime}_2}\ldots)$
is the associated partition of $n$ given by
$\alpha_1^\prime:=n-d-\sum_{i>0}\alpha_i$ and
$\alpha_i^\prime:=\alpha_{i-1}$ for $i\geq2$.
Here the assumption $n\geq 2d$ ensures that $\alpha'_1\geq 0$.

On the other hand, consider
\begin{equation}\label{eq:def_gamma_n}
\gamma_n^{\lambda}:=\Phi(\eps_n^\lambda)=
\Psi^{-1}\left(\prod_{i\geq1}c_i(pr_{1*}\ko_{\Xi_n})^{\alpha_i}\right).
\end{equation}
By the definition of the cup product, the elements
$\gamma_n^\lambda$ are all contained in $M$ and can therefore be
expressed as linear combinations of the $p_n^{\mu}$. We must show
that the associated coefficient matrix is invertible over $\IZ$.

This will be achieved by comparison with a third, rational basis of
the vector space $M\tensor_\IZ\IQ$, provided by the elements
\begin{equation}\label{eq:def_ch_n}
ch_n^{\lambda}:=
\Psi^{-1}\left(\prod_{i\geq1}ch_i(pr_{1*}\ko_{\Xi_n})^{\alpha_i}\right).
\end{equation}

\medbreak\noindent {\em Claim: Let $A$ be the matrix defined by
$ch_n^\lambda=\sum_{\mu\vdash d} A_{\mu\lambda}\gamma_n^\mu.$ Then}

\begin{equation}\label{eq:detA}
  \mid\det A\mid \ = \prod_{\lambda=(1^{\alpha_1}2^{\alpha_2}\ldots)\vdash d}\
\prod_{i\geq1}\left(\frac{1}{(i-1)!}\right)^{\alpha_i}
\end{equation}

Let $<$ be the order on the set of partitions
$\lambda=(1^{\alpha_1}2^{\alpha_2}\ldots)$ of $d$ corresponding to
the lexicographical order of the sequences
$(\alpha_1,\alpha_2,\ldots)$, so that for example the partition
$[1,1,\ldots,1]=(1^d)$ is the largest and $[d]=(d^1)$ is the
smallest.

Since Chern classes and the components of the Chern character
satisfy the universal identities
$$ch_k=\frac{(-1)^{k-1}}{(k-1)!}c_k+
\text{ polynomials in $c_1,\ldots,c_{k-1}$},$$
it follows that
$$ch_n^{\lambda}=\prod_{i\geq 1}\left(\frac{(-1)^{i-1}}{(i-1)!}
\right)^{\alpha_i}\, \gamma_n^\lambda+\text{ linear combination of
$\gamma_n^{\mu}$ with } \mu>\lambda.$$ This shows that $A$ is a lower
triangular matrix with diagonal entries
$$A_{\lambda\lambda}=\prod_{i\geq 1}\left(\frac{(-1)^{i-1}}{(i-1)!}
\right)^{\alpha_i}.$$
The claim follows directly from this.

\medbreak\noindent {\em Claim: Let $B$ be the matrix defined by
$ch_n^\lambda=\sum_{\mu\vdash d} B_{\mu\lambda} p_n^\mu$. Then
\begin{equation}\label{eq:detB}
\mid\det B\mid\ = \prod_{\lambda=(1^{\alpha_1}2^{\alpha_2}\ldots)\vdash
d} \ \prod_{i\geq 1}\left(\frac{1}{i!}\right)^{\alpha_i}\alpha_i!
\end{equation}
}

Recall that by Theorem \ref{th:BigDiffOperator} we have
$$\Psi^{-1}(ch_i(pr_{1*}\ko_{\Xi_n})\cup \Psi(y))=\kd_i(y)$$
for any polynomial $y\in \kp$, where the degree $i$ component
$\kd_i$ of the differential operator $\kd$ is given by
\begin{equation}\label{eq:DefOfD_i}
\kd_i=\frac{(-1)^{i}}{(i+1)!}\sum_{n_0,\ldots,n_i>0}p_{n_0+\ldots+n_i}
n_0\frac{\partial}{\partial p_{n_0}}\cdots n_i\frac{\partial}{\partial
p_{n_i}}.
\end{equation}
If $\kd_i$ is applied to a monomial $p_1^{\beta_1}\cdot\ldots\cdot
p_s^{\beta_s}$ with $\beta_1>i$, then the smallest component with
respect to the lexicographical order is that arising from the choice
$n_0=\ldots=n_i=1$ in (\ref{eq:DefOfD_i}). More precisely,
\begin{eqnarray*}\kd_i\left(\prod_{j\geq 1}\frac{1}{\beta_j!}
\left(\frac{p_j}{j}\right)^{\beta_j}\right)&=&
\frac{(-1)^i}{i!}(\beta_{i+1}+1)\frac{p_1^{\beta_1-i-1}}{(\beta_1-i-1)!}\\
&&\times\ \frac{1}{(\beta_{i+1}+1)!}
\left(\frac{p_{i+1}}{i+1}\right)^{\beta_{i+1}+1}\\
&&\times \prod_{j\neq 1,i+1}\frac{1}{\beta_j!}
\left(\frac{p_j}{j}\right)^{\beta_j} +\text{Terms of higher order}
\end{eqnarray*}
It follows by induction  that
\begin{eqnarray*}
ch_n^\lambda&=&\prod_{i\geq
  1}\kd_i^{\alpha_i}\left(\frac{p_1^n}{n!}\right)\\
&=&\prod_i\alpha_i!\left(\frac{(-1)^{\alpha_i}}{i!}\right)^{\alpha_i}
\cdot p^\lambda_n+\text{ linear combinations of $p_n^\mu$ with
$\mu^\prime>\lambda^\prime$}
\end{eqnarray*}
This shows that $B$ is a lower triangular matrix -- if we reorder
the $p_n^\lambda$ according to
$\mu\succ\lambda:\Leftrightarrow\mu^\prime>\lambda^\prime$ -- with
diagonal entries
$$B_{\lambda\lambda}=\prod_{i\geq 1}
\alpha_i!\left(\frac{(-1)^{\alpha_i}}{i!}\right)^{\alpha_i}.$$
The claim follows from this.

\medbreak\noindent
{\em Claim:
$$\left|\frac{\det A}{\det B}\right|=
\prod_{\lambda=(1^{\alpha_1}2^{\alpha_2}\ldots)\vdash d}\
\prod_{i\geq 1}\frac{i^{\alpha_i}}{\alpha_i!}=1.$$}

Of these two equalities the first is an immediate consequence of the two
previous claims. The second is a well known identity. In fact, it amounts
to realizing that each integer $k\in\{1,\ldots,d\}$ appears both in the
numerator and denominator with multiplicity
$$p(d-k)+p(d-2k)+p(d-3k)+\ldots,$$
where $p(s)$ is the number of partitions of $s$.
\end{proof}

\begin{remark}--- We have seen that for any
$\lambda=(1^{\alpha_1}2^{\alpha_2}\ldots)$ of $d$ there is a polynomial
$r_\lambda\in\IZ[c_1,c_2,\ldots]$ of (weighted) degree $d$ such that
$$\Psi\Phi(\chi_{\lambda'})=r_\lambda(c_1(pr_{1*}\ko_{\Xi_n}),\ldots,
c_d(pr_{1*}\ko_{\Xi_n}))
$$
whenever $n\geq 2d$. On the other hand, the kernel of the restriction map
$$\rho:\kc(S_{n+1})\to \kc(S_n)$$
is generated by all $\chi_{\lambda^\prime}$, where the coefficient
$\alpha_1^\prime$ in the presentation
$\lambda^\prime=(1^{\alpha^\prime_1}2^{\alpha^\prime_2}\ldots)$ vanishes.
This yields the following description of the cohomology ring in terms of
generators and relations:
$$
H^*(\Hilbn;\IZ)=\IZ[c_1,c_2,\ldots]/
(r_\lambda)_{\{\lambda\,|\,\sum_i(i+1)\alpha_i>n\}}.
$$
The polynomials $r_\lambda$ can be explicitly computed in $\kc(S_n)$ and have
a direct geometric interpretation: So for example among the first relations
that appear is $r_{(2^m)}$, where $m=\lceil (n+1)/2\rceil$, reflecting the
fact that the locus of points in $\Hilbn$ where more than $n/2$ pairs of
points collide is empty.
\end{remark}

%%% Biblio ---------------------------------------------------------------

\bibliographystyle{plain}

\parindent=0pt
\end{document}